\documentstyle[12pt]{article}
\def\Bbb R{{\rm \bf R}}
\def\proclaim#1{\vskip2mm{\bf #1}\em}
\def\endproclaim{\em \vskip2mm}
\def\tag#1{\eqno(#1)}
\def\gathered{\begin{array}{c}}
\def\endgathered{\end{array}}
\def\text{\mbox}

\begin{document}

\title {The Herglotz wave function, the Vekua transform\\
and the enclosure method}
\author{Masaru IKEHATA\\
Department of Mathematics,
Faculty of Engineering\\
Gunma University, Kiryu 376-8515, JAPAN}
\date{15 March 2005 Final}
\maketitle
\begin{abstract}
This paper gives applications of the enclosure method introduced by the author
to typical inverse obstacle and crack scattering problems in two dimensions.
Explicit extraction formulae of the convex hull of unknown
polygonal sound-hard obstacles and piecewise linear cracks from
the far field pattern of the scattered field at a fixed wave
number and at most two incident directions are given.
The main new points of this paper are: a combination of the enclosure method
and the Herglotz wave function; explicit construction of the density in the Herglotz wave function
by using the idea of the Vekua transform.
\noindent
By virtue of the construction,
one can avoid any restriction on the wave number in the extraction
formulae.  An attempt for the case when the far field pattern is
given on limited angles is also given.

\noindent
AMS: 35R30

\end{abstract}

\section{Introduction}

In this paper we consider typical inverse problems
for the Helmholtz equation in two-dimensions.
First let us consider an inverse obstacle scattering problem.
Let $D\subset\Bbb R^2$ be a bounded open set
with Lipschitz boundary and satisfy that $\Bbb
R^2\setminus\overline D$ is connected.
Using a variational method (e.g., see \cite{H} and references therein),
one knows that: given $k>0$ and $d\in S^1$ there exists a
unique $u\in C^{\infty}(\Bbb R^2\setminus\overline D)$ satisfying
(1) to (3) described below:

\noindent
(1) $u$ satisfies the Helmholtz equation
$$\displaystyle
\triangle u+k^2u=0\,\,\text{in}\,\Bbb R^2\setminus\overline D;
$$

\noindent
(2) there exists a disc $B_R$ with radius $R$ centered at $0$ such that
$\overline D\subset B_R$, $u\vert_{B_R\setminus\overline D}\in H^1(B_R\setminus\overline D)$
and for all $\phi\in H^1(B_R\setminus\overline D)$ with $\phi=0$ on $\partial B_R$
$$\displaystyle
\int_{B_R\setminus\overline D}(\nabla u\cdot\nabla\phi-k^2u\phi)dx=0;
\tag {1.1}
$$

\noindent
(3) $w=u-e^{ikx\cdot d}$ satisfies the outgoing Sommerfeld radiation condition
$$\displaystyle
\lim_{r\longrightarrow\infty}\sqrt{r}(\frac{\partial w}{\partial r}-ikw)=0
$$
where $r=\vert x\vert$.

\noindent
Note that the condition
$u\vert_{B_R\setminus\overline D}\in H^1(B_R\setminus\overline D)$ of (2) gives a restriction
on the singularity of $u$ in a neighbourhood of $\partial D$;
(1.1) is a weak formulation of the boundary condition $\partial u/\partial\nu=0$ on $\partial D$
where $\nu$ denotes the unit outward normal relative to $\Bbb R^2\setminus\overline D$.

It is well known that $w$ has the asymptotic expansion
as $r\longrightarrow\infty$ uniformly with respect to $\varphi\in\,S^1$:
$$
w(r\varphi)=\frac{e^{ikr}}{\sqrt{r}}F(\varphi;d,k)+O(\frac{1}{r^{3/2}}).
$$
The coefficient $F(\varphi;d,k)$ is called the far field pattern of
the reflected wave $w$ at direction $\varphi$.

In this paper we are interested to seek extraction formulae of information about the location and shape
of $D$ from the far field pattern for fixed $d$ and $k$.
Recently the author established an extraction formula of the convex hull of $D$
from the data $w$ on the boundary of any fixed open disc that contains $\overline D$
provided $D$ is polygonal (\cite{I3}). It is an application of the idea of the enclosure method
to the inverse scattering problem (see \cite{I2, IO} for a simpler problem).
Since the extraction formula in \cite{I3} is the starting point, we give a precise description
of the formula and point out the problem.

{\bf\noindent Definition 1.1.}
We say that $D$ is polygonal if $D=D_1\cup D_2\cup\cdots\cup D_m$; each $D_j$ is a simply
connected open set and polygon; $\overline D_j\cap\overline D_{j'}=\emptyset$ for $j\not=j'$.

In this paper we always assume that $D$ is polygonal.
Given direction $\omega\in S^1$ define $h_D(\omega)=\sup_{x\in\,D}x\cdot\omega$.
We call the function $\omega\mapsto h_D(\omega)$ the support function of $D$.  From the support function of $D$ one obtains the convex hull
of $D$.

{\bf\noindent Definition 1.2.}
We say that the direction $\omega$ is regular with respect to $D$ if the set
$\{x\,\vert x\cdot\omega=h_D(\omega)\}\cap\partial D$ consists of only one
point.

\noindent
We have already established the following.

\proclaim{\noindent Theorem 1.1(\cite{I3}).}
Let $\omega$ be regular with respect to $D$.  Then the formula
$$\displaystyle
\lim_{\tau\longrightarrow\infty}\frac{\displaystyle\log\vert
\int_{\partial B_R}(\frac{\partial u}{\partial \nu}v-\frac{\partial v}{\partial\nu}u)d\sigma\vert}
{\tau}
=h_D(\omega),
\tag {1.2}
$$
is valid where $v(x)=e^{x\cdot(\tau\omega+i\sqrt{\tau^2+k^2}\omega^{\perp})}$ and
$\omega^{\perp}=(\omega_2,-\omega_1)$ for $\omega=(\omega_1,\omega_2)$.  Moreover we have:

if $t\ge h_D(\omega)$, then
$$\displaystyle
\lim_{\tau\longrightarrow\infty}e^{-\tau t}
\vert
\int_{\partial B_R}(\frac{\partial u}{\partial \nu}v-\frac{\partial v}{\partial\nu}u)d\sigma\vert=0;
$$

if $t<h_D(\omega)$, then
$$\displaystyle
\lim_{\tau\longrightarrow\infty}e^{-\tau t}
\vert
\int_{\partial B_R}(\frac{\partial u}{\partial \nu}v-\frac{\partial v}{\partial\nu}u)d\sigma\vert=\infty.
$$

\endproclaim

It is well known that $\partial u/\partial\nu$ on $\partial B_R$
can be computed from $u$ on $\partial B_R$ by using the
Dirichlet-to-Neumann map outside $B_R$.  Moreover one can
calculate $u$ on $\partial B_R$ from $F(\,\cdot\,;d,k)$ for fixed
$d$ and $k$ by using, e.g., a formula in the point source method
(see \cite{P}). Thus one can say that (1.2) essentially gives an
extraction formula of the support function from the far field
pattern for fixed $d$ and $k$.  However, the computation involves
mainly two limiting procedures: the first is the procedure for
calculating $u$ and $\partial u/\partial\nu$ on $\partial B_R$
from $F(\,\cdot\,;d,k)$; the second is the procedure for
calculating $h_D(\omega)$ from $u$ and $\partial u/\partial\nu$ on
$\partial B_R$ by using (1.2).

In this paper we present a direct formula that extracts the value
of the support function of unknown polygonal sound hard obstacles
from the far field pattern for fixed $d$ and
$k$ just by using only one limiting procedure.

We identify the point $\vartheta=(\vartheta_1,\vartheta_2)\in S^1$
with the complex number $\vartheta_1+i\vartheta_2$ and denote it by the same
symbol $\vartheta$.

\noindent
Now we state the results.
\noindent
Given $N=1,\cdots$, $\tau>0$, $\omega\in S^1$ and $k>0$ define the function $g_N(\,\cdot\,;\tau,k,\omega)$ on $S^1$
by the formula
$$
\displaystyle
g_{N}(\varphi;\tau,k,\omega)
=\frac{1}{2\pi}
\sum_{\vert m\vert\le N}
\{\frac{ik\varphi}{(\tau+\sqrt{\tau^2+k^2})\omega}\}^m.
\tag {1.3}
$$

\proclaim{\noindent Theorem 1.2.}  Let $\omega$ be regular with
respect to $D$.  Let $\beta_0$ be the unique positive solution of the equation
$$\displaystyle
\frac{2}{e}s+\log s=0.
$$
Let $\beta$ satisfy $0<\beta<\beta_0$.  Let $\{\tau(N)\}_{N=1,\cdots}$ be an arbitrary
sequence of positive numbers satisfying, as $N\longrightarrow\infty$
$$\displaystyle
\tau(N)=\frac{\beta N}{eR}+O(1).
$$
Then the formula
$$\displaystyle
\lim_{N\longrightarrow\infty}
\frac{\displaystyle\log\vert\int_{S^1}F(-\varphi;d,k)g_N(\varphi;\tau(N),k,\omega)d\sigma(\varphi)\vert}
{\tau(N)}=h_D(\omega),
\tag {1.4}
$$
is valid.
Moreover we have:

if $t\ge h_D(\omega)$, then
$$\displaystyle
\lim_{N\longrightarrow\infty}e^{-\tau(N)t}
\vert\int_{S^1}F(-\varphi;d,k)g_N(\varphi;\tau(N),k\,\omega)d\sigma(\varphi)\vert=0;
$$

if $t<h_D(\omega)$, then
$$\displaystyle
\lim_{N\longrightarrow\infty}e^{-\tau(N)t}
\vert\int_{S^1}F(-\varphi;d,k)g_N(\varphi;\tau(N),k,\omega)d\sigma(\varphi)\vert=\infty.
$$

\endproclaim

\noindent Note that there is no restriction on $k$.
(1.4) means that if one knows the Fourier coefficients of the far field pattern
$$\displaystyle
\int_{S^1}F(\varphi;d,k)\varphi^md\sigma(\varphi)=(-1)^m\int_{S^1}F(-\varphi;d,k)\varphi^md\sigma(\varphi)
, \,\,\vert m\vert\le N
$$
for sufficiently large $N$ and a disc that contains $\overline D$, then one can know an approximate value of $h_D(\omega)$.

\noindent
$g=g_N(\,\cdot\,;\tau(N),k,\omega)$ satisfies, as
$N\longrightarrow\infty$
$$\displaystyle
\int_{S^1}e^{iky\cdot\varphi}g(\varphi)d\sigma(\varphi)
\approx e^{y\cdot(\tau(N)\omega+i\sqrt{\tau(N)^2+k^2}\omega^{\perp})}, \,\,y\in\overline B_R
$$
in an appropriate sense (see Theorem 2.1 in Section 2).
The function
$$
\int_{S^1}e^{iky\cdot\varphi}g(\varphi)d\sigma(\varphi)
$$
is called the Herglotz wave function with density $g$ and
satisfies the Helmholtz equation $\triangle v+k^2 v=0$ in the
whole space.

For a simple explanation of the origin of the desired density we make use of the idea of
the Vekua transform which maps harmonic functions into solutions of
the Helmholtz equations (see \cite{V, BG} and also \cite{CD} for a recent study).
Recall the Bessel function of order $m=0,\pm 1,\pm 2,\cdots$ given by the formula
$$\displaystyle
J_m(z)=
(\frac{z}{2})^m\sum_{n=0}^{\infty}
\frac{(-1)^n}
{n!\Gamma(n+1+m)}(\frac{z}{2})^{2n}.
$$
The Vekua transform in two-dimensions takes the form:
$$\displaystyle
T_kv(x)=v(x)-\frac{k\vert x\vert}{2}\int_0^1 v(tx)J_1(k\vert x\vert\sqrt{1-t})\frac{dt}{\sqrt{1-t}}
$$
where $v$ is an arbitrary harmonic function in $\Bbb R^2$.  Using the formulae
$$
\displaystyle
k\vert x\vert\int_0^1(1-w^2)^mJ_1(k\vert x\vert w)dw
=1-(\frac{2}{k\vert x\vert})^m m!J_m(k\vert x\vert),\,\,m=0,1,\cdots
$$
one can easily know that, for $m=0,1,2,\cdots$
$$
\displaystyle
T_k: r^m e^{\pm im\theta}\mapsto (\frac{2}{k})^m m!J_m(kr)e^{\pm im\theta}.
$$
Using this property (rule), we reduce the construction problem of the density to
that of the density in the "harmonic" Herglotz wave function
$$\displaystyle
\int_{S^1}\{e^{ik\overline\varphi(y_1+iy_2)/2}
+e^{ik\varphi(y_1-iy_2)/2}-1\}g(\varphi)d\sigma(\varphi).
$$
See Section 2 for details.  Note that one can state and prove the
main result of Section 2 (Theorem 2.1) without using the Vekua
transform.  However, the presence of the "harmonic" Herglotz wave
function is something interesting.

\noindent

\noindent
(1.4) gives us a direct formula for extracting information about the convex hull from
the far field pattern by using only one limiting procedure.
An interesting question is to seek such a formula in the case when $F(\,\cdot\,;d,k)$ is given only on
a proper subset of $S^1$. In Section 4 we present a small first step
that gives extraction formulae of the support function from
the far field pattern on given nonempty open subset of $S^1$ provided the center
of the coordinates is inside the obstacles.

For another approach to the problem one can cite the no response
test introduced in \cite{LP}.  Therein, for fixed $\epsilon>0$
they define the functional of test domains $D_t$ by taking the least
upper bound of the quantity
$$\displaystyle
\vert\int_{-\Gamma}F(-\varphi;d,k)g(\varphi)d\sigma(\varphi)\vert
$$
with respect to all $g\in L^2(-\Gamma)$ satisfying
$$
\displaystyle
\sup_{y\in
D_t}
\vert\int_{-\Gamma}e^{iky\cdot\varphi}g(\varphi)d\sigma(\varphi)\vert\le\epsilon.
\tag {1.5}
$$
Their idea is to make use of this functional to decide whether
$D$ is contained in $D_t$. In contrast to the enclosure method the
form of the density of the Herglotz wave function is not specified
except for a restriction on the bound (1.5) of the corresponding
Herglotz wave function. Note that Theorem 1.2 gives an explicit
way in the case when $\Gamma=S^1$ of deciding whether $D$ is
contained in the special test domain $B_R\cap\{x\in\Bbb
R^2\,\vert\,x\cdot\omega<t\}$ for each $t$ with $\vert t\vert<R$.

\noindent
Note that if one has the complete knowledge of the far
field pattern for all $d$ and fixed $k$, then one has two
reconstruction formulae for general $D$. One is due to Kirsch
(\cite{Ki}) and another is due to the author (\cite{I1, I11}).

Next we consider an inverse scattering problem for piecewise linear cracks.
Let $\Sigma$ be the union of finitely many disjoint closed piecewise linear segments
$\Sigma_1, \Sigma_2,\cdots,\Sigma_m$.  Assume that there exists  a simply connected open set $D$
such that $D$ is a polygon and each $\Sigma_j$ is contained in $\partial D$.
We assume that $\overline D\subset B_R$.
We denote by $\nu$ the unit outward normal relative to $B_R\setminus\overline D$.
Denote by $H^1(B_R\setminus\Sigma)$
the set of all $L^2(B_R)$ functions $u$ such that, $u^+=u\vert_{B_R\setminus\overline D}
\in H^1(B_R\setminus\overline D)$, $u^-=u\vert_D\in H^1(D)$ and
$u^+=u^-$ on $\partial D\setminus\Sigma$.

\noindent
It is well known that: given $k>0$ and $d\in S^1$ there exists a unique $u\in C^{\infty}(\Bbb R^2\setminus\Sigma)$
satisfying (4) to (6) described below:

\noindent
(4)  $u$ satisfies the Helmholtz equation
$$
\displaystyle
\triangle u+k^2u=0\,\,\text{in $\Bbb R^2\setminus\Sigma$;}
$$

\noindent
(5)  $u\vert_{B_R}\in H^1(B_R\setminus\Sigma)$ and for all $\phi\in H^1(B_R\setminus\Sigma)$
with $\phi=0$ on $\partial B_R$
$$
\displaystyle
\int_{B_R\setminus\Sigma}(\nabla u\cdot\nabla\varphi-k^2u\varphi)dx=0;
\tag {1.6}
$$

\noindent
(6)  $w=u-e^{ikx\cdot d}$ satisfies the outgoing Sommerfeld radiation condition (see (3)).

\noindent
(1.6) is a weak formulation of the boundary condition $\partial u/\partial\nu=0$
on $\Sigma$ where $\nu$ stands for the unit outward normal
relative to $B_R\setminus\overline D$.

Define
$$\displaystyle
h_{\Sigma}(\omega)=\sup_{x\in \Sigma}x\cdot\omega,\,\,\omega\in\,S^1.
$$
We call $h_{\Sigma}$ the support function of $\Sigma$.
We say that the direction $\omega\in S^1$ is
regular with respect to $\Sigma$ if the set
$\{x\,\vert\,x\cdot\omega=h_{\Sigma}(\omega)\}\cap\Sigma$ consists of
only one point.
\noindent

Given $\omega=(\omega_1,\omega_2)\in S^1$ set
$\omega^{\perp}=(\omega_2,-\omega_1)$.
Define the indicator function $I_{\omega}^d(\tau,t)$ by the formula
$$
I^d_{\omega}(\tau,t)
=e^{-\tau t}\vert\int_{\vert x\vert=R}
(\frac{\partial u}{\partial\nu}v-\frac{\partial v}{\partial\nu}u)d\sigma\vert
$$
where $-\infty<t<\infty$, $\tau>0$ and $v=e^{x\cdot(\tau\omega+i\sqrt{\tau^2+k^2}\omega^{\perp})}$.

In \cite{I3} we gave an extraction formula of the convex hull of
$\Sigma$ from the indicator function for a single incident
direction $d$.

\proclaim{\noindent Theorem 1.3(\cite{I3}).} Let $\omega$ be regular with respect to $\Sigma$.
If every end point of $\Sigma_1,\Sigma_2,\cdots,\Sigma_m$ satisfies
$x\cdot\omega<h_{\Sigma}(\omega)$, then the formula
$$\displaystyle
\lim_{\tau\longrightarrow\infty}\frac{\log I^d_{\omega}(\tau,0)}
{\tau}=h_{\Sigma}(\omega),
$$
is valid. Moreover, we have:

if $t\ge h_{\Sigma}(\omega)$, then $\lim_{\tau\longrightarrow\infty}I^d_{\omega}(\tau,t)=0$;

if $t<h_{\Sigma}(\omega)$, then $\lim_{\tau\longrightarrow\infty}I^d_{\omega}(\tau,t)=\infty$.

\noindent
If there is an end point $y$ of some $\Sigma_j$ such that
$y\cdot\omega=h_{\Sigma}(\omega)$, then, for $d$ that is not
perpendicular to $\nu$ on $\Sigma_j$ near the point, the same
conclusions as above are valid.

\endproclaim
\noindent
Note that $\nu$ on $\Sigma_j$ in a neighbourhood of $y$
becomes a constant vector if $y$ is an end point of $\Sigma_j$.

\noindent
From the definition of regularity of $\omega$ one knows that the set
$\{x\,\vert\,x\cdot\omega=h_{\Sigma}(\omega)\}\cap\Sigma$ consists of
a single point (say $x_0$).  Then every endpoint of $\Sigma_1,\cdots,\Sigma_m$
satisfies $x\cdot\omega<h_{\Sigma}(\omega)$ if and only if there exists $j$ such that $x_0\in\Sigma_j$
and $\nu$ on $\Sigma_j$ has discontinuity at $x_0$.

One of two implications of Theorem 1.3 is:  one can delete the sentences "
If every end point of $\Sigma_1,\Sigma_2,\cdots,\Sigma_m$ satisfies
$x\cdot\omega<h_{\Sigma}(\omega)$" and "If there is an end point...are valid" in Theorem 1.3
by introducing the new indicator function
$$
I_{\omega}(\tau,t)=I_{\omega}^{d_1}(\tau,t)+I_{\omega}^{d_2}(\tau,t)
$$
where $d_1, d_2$ are two arbitrary linearly independent incident directions and fixed.  
This indicator function makes use of
two reflected waves for two incident plane waves at a fixed wave number.   This idea is coming from that of
the multi probe method (see \cite{IP}).
We obtain the following.

\proclaim{\noindent Theorem 1.4.} Let $\omega$ be regular with respect to $\Sigma$.
The formula
$$\displaystyle
\lim_{\tau\longrightarrow\infty}\frac{\log I_{\omega}(\tau,0)}
{\tau}=h_{\Sigma}(\omega),
$$
is valid. Moreover, we have:

if $t\ge h_{\Sigma}(\omega)$, then $\lim_{\tau\longrightarrow\infty}I_{\omega}(\tau,t)=0$;

if $t<h_{\Sigma}(\omega)$, then $\lim_{\tau\longrightarrow\infty}I_{\omega}(\tau,t)=\infty$.

\endproclaim

\noindent
The next theorem corresponds to Theorem 1.2.

\proclaim{\noindent Theorem 1.5.} Fix two arbitrary linearly independent directions $d_1$ and $d_2$.
Let $\omega$ be regular with
respect to $\Sigma$.  Let $\{\tau(N)\}_{N=1,\cdots}$ be the same as that of Theorem 1.2.
Then the formula
$$\displaystyle
\lim_{N\longrightarrow\infty}
\frac{\displaystyle\log\{
\sum_{j=1}^2
\vert\int_{S^1}F(-\varphi;d_j,k)g_N(\varphi;\tau(N),k,\omega)d\sigma(\varphi)\vert\}}
{\tau(N)}=h_{\Sigma}(\omega),
$$
is valid.  Moreover we have:

if $t\ge h_{\Sigma}(\omega)$, then
$$\displaystyle
\lim_{N\longrightarrow\infty}e^{-\tau(N)t}
\sum_{j=1}^2\vert\int_{S^1}F(-\varphi;d_j,k)g_N(\varphi;\tau(N),k,\omega)d\sigma(\varphi)\vert=0;
$$

if $t<h_{\Sigma}(\omega)$, then
$$\displaystyle
\lim_{N\longrightarrow\infty}e^{-\tau(N)t}
\sum_{j=1}^2\vert\int_{S^1}F(-\varphi;d_j,k)g_N(\varphi;\tau(N),k,\omega)d\sigma(\varphi)\vert=\infty.
$$

\endproclaim

The outline of this paper is as follows.  In Section 2 we describe
the derivation of the density given by (1.3) by using the idea of
the Vekua transform and establish the desired property. The proof
of Theorem 1.2, 1.4 and comment on the proof of Theorem 1.5 are
found in Section 3.  In Section 4 a case when the far field
pattern is given on limited angles is discussed.

\section{Construction of the density in the Herglotz wave function}

The aim of this section is to construct a density $g\in L^2(S^1)$
explicitly such that
$$\displaystyle
\int_{S^1}e^{iky\cdot\varphi}g(\varphi)d\sigma(\varphi)
\approx e^{y\cdot(\tau\omega+i\sqrt{\tau^2+k^2}\omega^{\perp})}, \,\,y\in\overline B_R.
$$

\noindent
We start with an elementary fact.

\proclaim{\noindent Proposition 2.1.} Let $m$ be an arbitrary
integer.  Let $\xi$ be an arbitrary complex vector. We have
$$
\frac{1}{2\pi}\int_{S^1}e^{\vartheta\cdot\xi}(\vartheta_1+i\vartheta_2)^md\sigma(\vartheta)
=\left\{
\begin{array}{lr}
\displaystyle
(\frac{z}{2})^{m}
\sum_{n=0}^{\infty}
\frac{1}{n!(m+n)!}
(\frac{1}{2})^{2n}(z^*z)^n, & \quad\text{if $m\ge 0$,}\\
\\
\displaystyle
(\frac{z^*}{2})^{-m}
\sum_{n=0}^{\infty}
\frac{1}{n!(-m+n)!}
(\frac{1}{2})^{2n}(z^*z)^n, & \quad\text{if $m<0$}
\end{array}
\right.
\tag {2.1}
$$
where $\xi=(\xi_1,\xi_2)$, $z=\xi_1+i\xi_2$ and $z^*=\xi_1-i\xi_2$.
\endproclaim

{\it\noindent Proof.}
Let $w=e^{i\theta}$.
Then
$$\begin{array}{c}
\displaystyle
\xi_1\cos\theta+\xi_2\sin\theta
=\frac{\xi_1}{2}(w+w^{-1})+\frac{\xi_2}{2i}(w-w^{-1})\\
\\
\displaystyle
=\frac{1}{2}(z^*w+zw^{-1}).
\end{array}
$$
This gives the expression
$$\begin{array}{c}
\displaystyle
\int_{S^1}e^{\vartheta\cdot\xi}(\vartheta_1+i\vartheta_2)^md\sigma(\vartheta)
=\int_0^{2\pi}e^{\xi_1\cos\theta+\xi_2\sin\theta}e^{im\theta}d\theta\\
\\
\displaystyle
=\frac{1}{i}\int_{\vert w\vert=1}w^{m-1}e^{\frac{1}{2}(z^*w+zw^{-1})}dw\\
\\
\displaystyle
=2\pi\text{Res}_{w=0}(w^{m-1}e^{\frac{1}{2}(z^*w+zw^{-1})}).
\end{array}
$$
Write
$$\begin{array}{c}
\displaystyle
w^{m-1}e^{\frac{1}{2}(z^*w+zw^{-1})}
=\sum_{n=0}^{\infty}\frac{1}{n!}(\frac{1}{2})^n
(z^*w+zw^{-1})^nw^{m-1}\\
\\
=\displaystyle
\sum_{n=0}^{\infty}\frac{1}{n!}(\frac{1}{2})^n\sum_{r=0}^n
\left(
\begin{array}{c}
 n \\r
\end{array}
\right)
(z^*)^{n-r}z^r w^{n-2r+m-1}.
\end{array}
\tag {2.2}
$$
Consider the case when $m\ge 0$.
If $n<m$, then $n+m>2n\ge 2r$ for $0\le r\le n$.  Thus $n-2r+m-1\not=-1$.
This gives
$$\begin{array}{c}
\displaystyle
\text{Res}_{w=0}
(w^{m-1}e^{\frac{1}{2}(z^*w+zw^{-1})})
=\text{Res}_{w=0}
(\sum_{n=m}^{\infty}\frac{1}{n!}(\frac{1}{2})^n
\sum_{r=0}^n
\left(
\begin{array}{c}
n\\r
\end{array}
\right)
(z^*)^{n-r}z^rw^{n+m-2r-1})\\
\\
\displaystyle
=\text{Res}_{w=0}
(\sum_{l=0}^{\infty}
\frac{1}{(m+l)!}
(\frac{1}{2})^{m+l}
\sum_{r=0}^{m+l}
\left(
\begin{array}{c}
m+l\\ r
\end{array}
\right)
(z^*)^{m+l-r}z^rw^{2(m-r)+l-1}).
\end{array}
$$
If $l=2l', l'=0,\cdots$, then
$2(m-r)+l-1=2(m+l'-r)-1$.  Since $0\le r\le m+2l'$,
$2(m-r)+l-1$ becomes $-1$ only for $r=m+l'$.
If $l=2l'+1, l'=0,\cdots$, then
$2(m-r)+l-1=2(m+l'-r)$ and thus never become $-1$.
Therefore
$$\begin{array}{c}
\displaystyle
\text{Res}_{w=0}
(w^{m-1}e^{\frac{1}{2}(z^*w+zw^{-1})})
=\sum_{l'=0}^{\infty}
\frac{1}{(m+2l')!}
(\frac{1}{2})^{m+2l'}
\left(
\begin{array}{c}
m+2l'\\m+l'\end{array}
\right)
(z^*)^{m+2l'-(m+l')}z^{m+l'}\\
\\
=\displaystyle
\sum_{l'=0}^{\infty}
\frac{1}{l'!(m+l')!}
(\frac{1}{2})^{m+2l'}
(z^*)^{l'}z^{m+l'}\\
\\
\displaystyle
=(\frac{z}{2})^m
\sum_{n=0}^{\infty}\frac{1}{n!(m+n!)}(\frac{1}{2})^{2n}
(z^*z)^n.
\end{array}
\tag {2.3}
$$
This gives (2.1) in the case when $m\ge 0$.
Consider the case when $m<0$.  If $0\le n<-m$, then $n+m<0$ and thus
$n+m-1-2r<-2r-1\le -1$ for $0\le r\le n$.  Then from (2.2) we have
$$\begin{array}{c}
\displaystyle
\text{Res}_{w=0}
(w^{m-1}e^{\frac{1}{2}(z^*w+zw^{-1})})
=\text{Res}_{w=0}
(\sum_{n=-m}^{\infty}
\frac{1}{n!}(\frac{1}{2})^n
\sum_{r=0}^n\left(\begin{array}{c}
n\\ r\end{array}\right)
(z^*)^{n-r}z^rw^{n+m-2r-1})\\
\\
\displaystyle
=\text{Res}_{w=0}
(\sum_{l=0}^{\infty}
\frac{1}{(-m+l)!}
(\frac{1}{2})^{-m+l}
\sum_{r=0}^{-m+l}
\left(\begin{array}{c}
-m+l\\ r\end{array}
\right)(z^*)^{-m+l-r}z^rw^{l-2r-1})\\
\\
\displaystyle
=\text{Res}_{w=0}
(\sum_{l'=0}^{\infty}
\frac{1}{(-m+2l')!}(\frac{1}{2})^{-m+2l'}
\sum_{r=0}^{-m+2l'}
\left(\begin{array}{c}
-m+2l'\\ r\end{array}
\right)(z^*)^{-m+2l'-r}z^r
w^{2(l'-r)-1})\\
\\
\displaystyle
=\sum_{l'=0}^{\infty}
\frac{1}{(-m+2l')!}
(\frac{1}{2})^{-m+2l'}
\left(\begin{array}{c}
-m+2l'\\ l'\end{array}
\right)(z^*)^{-m+l'}z^{l'}\\
\\
\displaystyle
=\sum_{l'=0}^{\infty}\frac{1}{l'!(-m+l')!}
(\frac{1}{2})^{-m+2l'}
z^{l'}(z^*)^{-m+l'}\\
\\
\displaystyle
=(\frac{z^*}{2})^{-m}
\sum_{l'=0}^{\infty}\frac{1}{l'!(-m+l'!)}(\frac{1}{2})^{2l'}(z^*z)^{l'}.
\end{array}
$$

\noindent
$\Box$

\noindent
Given $\omega=(\omega_1,\omega_2)\in S^1$ set $\omega^{\perp}=(\omega_2,-\omega_1)$.
In this paper we denote the complex number $\omega_1+i\omega_2$ by $\omega$ again and thus
$\overline\omega=\omega_1-i\omega_2$.

We give another expression of (2.1) for $\xi=r(\tau\omega+i\sqrt{\tau^2+k^2}\omega^{\perp})$ with $r>0$.
Since $\xi_1=r(\tau\omega_1+i\sqrt{\tau^2+k^2}\omega_2)$ and $\xi_2=r(\tau\omega_2-i\sqrt{\tau^2+k^2}\omega_1)$,
we have
$$\displaystyle
z=r(\tau+\sqrt{\tau^2+k^2})\omega,\,\,
z^*=r(\tau-\sqrt{\tau^2+k^2})\overline\omega
$$
and $z^*z=-(rk)^2$.  Then (2.1) gives
$$
\frac{1}{2\pi}\int_{S^1}e^{\vartheta\cdot\xi}(\vartheta_1+i\vartheta_2)^md\sigma(\vartheta)
=\left\{
\begin{array}{lr}
\displaystyle
\{\frac{(\tau+\sqrt{\tau^2+k^2})\omega}{k}\}^mJ_m(kr), & \quad\text{if $m\ge 0$,}\\
\\
\displaystyle
\{\frac{(\tau-\sqrt{\tau^2+k^2})\overline\omega}{k}\}^{-m}J_{-m}(kr), & \quad\text{if $m<0$}
\end{array}
\right.
\tag {2.4}
$$
Note that using the equation $J_{-m}(kr)=(-1)^mJ_m(kr)$ one can rewrite (2.4) as
for all integer $m$
$$\displaystyle
\frac{1}{2\pi}\int_{S^1}e^{r\vartheta\cdot(\tau\omega+i\sqrt{\tau^2+k^2}\omega^{\perp})}
(\vartheta_1+i\vartheta_2)^m d\sigma(\vartheta)
=\{\frac{(\tau+\sqrt{\tau^2+k^2})\omega}{k}\}^mJ_m(kr),
$$
however, for our purpose (2.4) is more convenient.
From (2.4) we obtain the expansion formula
$$\begin{array}{c}
\displaystyle
e^{y\cdot(\tau\omega+i\sqrt{\tau^2+k^2}\omega^{\perp})}
=\sum_{m=0}^{\infty}
\{\frac{(\tau-\sqrt{\tau^2+k^2})\overline\omega}{k}\}^mJ_m(kr)e^{im\theta}\\
\\
\displaystyle
+\sum_{m=0}^{\infty}\{\frac{(\tau+\sqrt{\tau^2+k^2})\omega}{k}\}^mJ_m(kr)e^{-im\theta}
-J_0(kr)
\end{array}
\tag {2.5}
$$
where $y=(r\cos\,\theta,r\sin\,\theta)$.

\noindent
Define the harmonic function $e_{\omega}(y;\tau,k)$ in the whole space by the formula
$$
\displaystyle
e_{\omega}(y;\tau,k)
=e^{\displaystyle (\tau-\sqrt{\tau^2+k^2})\overline\omega(y_1+iy_2)/2}
+e^{\displaystyle (\tau+\sqrt{\tau^2+k^2})\omega(y_1-iy_2)/2}-1.
$$
Then from (2.5) one immediately obtains the following.

\proclaim{\noindent Proposition 2.2.} The Vekua transform of the
harmonic function $\displaystyle e_{\omega}(y;\tau,k)$ coincides
with $\displaystyle
e^{y\cdot(\tau\omega+i\sqrt{\tau^2+k^2}\omega^{\perp})}$.

\endproclaim

Next consider the case when $\xi$ in (2.1) is given by the complex vector
$\displaystyle \xi=ikr\varphi,\,\varphi\in S^1,\,\,r>0$.
Then $z=ikr\varphi$, $z^*=ikr\overline\varphi$ and $z^*z=-r^2k^2$.  Thus from (2.1) we have
$$
\frac{1}{2\pi}\int_{S^1}e^{ikr\vartheta\cdot\varphi}(\vartheta_1+i\vartheta_2)^md\sigma(\vartheta)
=\left\{
\begin{array}{lr}
\displaystyle
(i\varphi)^mJ_m(kr), & \quad\text{if $m\ge 0$,}\\
\\
\displaystyle
(i\overline\varphi)^{-m}J_{-m}(kr), & \quad\text{if $m<0$.}
\end{array}
\right.
$$
This gives the Jacobi-Anger expansion
$$
\displaystyle
e^{iky\cdot\varphi}
=\sum_{m=0}^{\infty}(i\overline\varphi)^mJ_{m}(kr)e^{im\theta}+
\sum_{m=0}^{\infty}(i\varphi)^mJ_m(kr)e^{-im\theta}-J_0(kr)
\tag {2.6}
$$
where $y=(r\cos\,\theta,r\sin\,\theta)$. Then the following
statement becomes trivial.

\proclaim{\noindent Proposition 2.3.}
The Vekua transform of the harmonic function
$$\displaystyle
e^{\displaystyle ik\overline\varphi(y_1+iy_2)/2}+e^{\displaystyle ik\varphi(y_1-iy_2)/2}-1
$$
coincides with $e^{iky\cdot\varphi}$.

\endproclaim

Let $\Gamma$ be a non empty open subset of $S^1$.
Given $g\in L^2(S^1)$ the function
$$\displaystyle
\int_{\Gamma}\{e^{\displaystyle ik\overline\varphi(y_1+iy_2)/2}+e^{\displaystyle ik\varphi(y_1-iy_2)/2}-1\}g(\varphi)d\sigma(\varphi)
$$
is harmonic in the whole space.  As a corollary of Proposition 2.3 one knows that
the Vekua transform of this harmonic function coincides with the Herglotz wave function with density $g$
$$\displaystyle
\int_{\Gamma}e^{iky\cdot\varphi}g(\varphi)d\sigma(\varphi).
$$

\noindent
Taking account of Proposition 2.2 and the fact mentioned above, it suffices to construct $g$ in such a way that
$$\begin{array}{c}
\displaystyle
\int_{\Gamma}\{e^{\displaystyle ik\overline\varphi(y_1+iy_2)/2}+e^{\displaystyle ik\varphi(y_1-iy_2)/2}-1\}g(\varphi)d\sigma(\varphi)\\
\\
\displaystyle
\approx
e_{\omega}(y;\tau,k).
\end{array}
\tag {2.7}
$$

\noindent
Using the power series expansion of the exponential function, one knows that if $g$ satisfies the system of equations
$$
\displaystyle
(\frac{ik}{2})^m\int_{\Gamma}(\overline\varphi)^mg(\varphi)d\sigma(\varphi)
=\{\frac{(\tau-\sqrt{\tau^2+k^2})\overline\omega}{2}\}^m, m=0, 1,\cdots
\tag {2.8}
$$
and
$$
\displaystyle
(\frac{ik}{2})^m\int_{\Gamma}\varphi^mg(\varphi)d\sigma(\varphi)
=\{\frac{(\tau+\sqrt{\tau^2+k^2})\omega}{2}\}^m, m=1,\cdots,
\tag {2.9}
$$
then $g$ satisfies (2.7) exactly.
We construct $g$ in the form
$$
g(\varphi)=\sum_{m=0}^{\infty}\beta_m\varphi^m+\sum_{m=1}^{\infty}\beta_{-m}\overline\varphi^{m}.
$$
Now consider the case when $\Gamma=S^1$.  Since
$$
\displaystyle
\frac{1}{2\pi}\int_{S^1}\overline\varphi^mg(\varphi)d\sigma(\varphi)=\beta_m
$$
and
$$
\displaystyle
\frac{1}{2\pi}\int_{S^1}\varphi^mg(\varphi)d\sigma(\varphi)
=\beta_{-m},
$$
from (2.8) and (2.9) we get
$$\displaystyle
\beta_m=\frac{1}{2\pi}\{\frac{ik}{(\tau+\sqrt{\tau^2+k^2})\omega}\}^m,\,\,m=0,\cdots
$$
and
$$
\displaystyle
\beta_{-m}
=\frac{1}{2\pi}\{\frac{(\tau+\sqrt{\tau^2+k^2})\omega}{ik}\}^m, m=1,\cdots.
$$
Then $g$ becomes
$$
\displaystyle
g(\varphi)
=\sum_{m=0}^{\infty}\frac{1}{2\pi}\{\frac{ik\varphi}{(\tau+\sqrt{\tau^2+k^2})\omega}\}^m
+\sum_{m=1}^{\infty}\frac{1}{2\pi}\{\frac{(\tau+\sqrt{\tau^2+k^2})\omega}{ik\varphi}\}^m
\tag {2.10}
$$
The first term is convergent and has the form
$$\displaystyle
\frac{1}{2\pi}\frac{(\tau+\sqrt{\tau^2+k^2})\omega}{(\tau+\sqrt{\tau^2+k^2})\omega-ik\varphi}.
$$
However the second term is always divergent since $\tau+\sqrt{\tau^2+k^2}>k$.  So we consider a
truncation of (2.10):
$$
\displaystyle
g_{N}(\varphi;\tau,k,\omega)
=\frac{1}{2\pi}\sum_{m=0}^N
\{\frac{ik\varphi}{(\tau+\sqrt{\tau^2+k^2})\omega}\}^m
+\frac{1}{2\pi}\sum_{m=1}^{N}\{\frac{(\tau+\sqrt{\tau^2+k^2})\omega}{ik\varphi}\}^m
$$
where $N=1,\cdots$.  This coincides with the expression given by (1.3).
Then one obtains
$$\begin{array}{c}
\displaystyle
\int_{S^1}\{e^{\displaystyle ik\overline\varphi(y_1+iy_2)/2}+e^{\displaystyle ik\varphi(y_1-iy_2)/2}-1\}g_N(\varphi;\tau,k,\omega)d\sigma(\varphi)
-e_{\omega}(y;\tau,k)\\
\\
\displaystyle
=-\sum_{m>N}\frac{1}{m!}
\{\frac{(\tau-\sqrt{\tau^2+k^2})\overline\omega}{2}\}^m(y_1+iy_2)^m
-\sum_{m>N}\frac{1}{m!}
\{\frac{(\tau+\sqrt{\tau^2+k^2})\omega}{2}\}^m(y_1-iy_2)^m.
\end{array}
\tag {2.11}
$$
This shows $g_{N}(\,\cdot\,;\tau,k,\omega)$ satisfies (2.7) in
this sense.  Taking the Vekua transform of the both sides of
(2.11) we obtain the equation
$$\begin{array}{c}
\displaystyle
\int_{S^1}e^{iky\cdot\varphi}g_N(\varphi;\tau,k,\omega)d\sigma(\varphi)
-e^{y\cdot(\tau\omega+i\sqrt{\tau^2+k^2}\omega^{\perp})}\\
\\
\displaystyle
=-\sum_{m>N}\{\frac{(\tau-\sqrt{\tau^2+k^2})\overline\omega}{k}\}^mJ_m(kr)e^{im\theta}
-\sum_{m>N}\{\frac{(\tau+\sqrt{\tau^2+k^2})\omega}{k}\}^mJ_m(kr)e^{-im\theta}
\end{array}
\tag {2.12}
$$
where $y=(r\cos\,\theta,r\sin\,\theta)$.  Note that this can be checked also directly.
For our purpose we have to consider how to choose $\tau$ depending on $N$.
One answer to this question is the following and it is the main result of this section.

\proclaim{\noindent Theorem 2.1.}
Let $\beta_0$ be the unique positive solution of the equation
$$\displaystyle
\frac{2}{e}s+\log s=0.
$$
Let $\beta$ satisfy $0<\beta<\beta_0$.  Let $\{\tau(N)\}_{N=1,\cdots}$ be an arbitrary
sequence of positive numbers satisfying, as $N\longrightarrow\infty$
$$\displaystyle
\tau(N)=\frac{\beta N}{eR}+O(1).
$$
Then we have, as $N\longrightarrow\infty$
$$\begin{array}{c}
\displaystyle
e^{R\tau(N)}\sup_{\vert y\vert\le R}\vert\int_{S^1}e^{iky\cdot\varphi}g_N(\varphi;\tau(N),k,\omega)d\sigma(\varphi)
-e^{y\cdot(\tau(N)\omega+i\sqrt{\tau(N)^2+k^2}\omega^{\perp})}\vert
\\
\\
\displaystyle
+e^{R\tau(N)}\sup_{\vert y\vert\le R}\vert\nabla
\{\int_{S^1}e^{iky\cdot\varphi}g_N(\varphi;\tau(N),k,\omega)d\sigma(\varphi)
-e^{y\cdot(\tau(N)\omega+i\sqrt{\tau(N)^2+k^2}\omega^{\perp})}\}\vert
\\
\\
\displaystyle
=O(N^{-\infty}).
\end{array}
\tag {2.13}
$$
\endproclaim

{\it\noindent Proof.} We give a direct proof without referring the
property of the Vekua transform. Set
$$
\begin{array}{c}
\displaystyle
R_N(y;\tau)=\sum_{m>N}\{\frac{(\tau-\sqrt{\tau^2+k^2})\overline\omega}{k}\}^mJ_m(kr)e^{im\theta}
\\
\\
\displaystyle
S_N(y;\tau)=\sum_{m>N}\{\frac{(\tau+\sqrt{\tau^2+k^2})\omega}{k}\}^mJ_m(kr)e^{-im\theta}
\\
\\
\displaystyle
E(\tau;N)=\frac{1}{N!}
\{\frac{R(\tau+\sqrt{\tau^2+k^2})}{2}\}^N
e^{R(\tau+\sqrt{\tau^2+k^2})/2}.
\end{array}
$$
The estimate
$$
\vert J_{m}(kr)\vert\le (\frac{kr}{2})^m\frac{1}{m!},
$$
is well known.  Then we have, for all $y$ with $\vert y\vert\le R$
$$
\displaystyle
\vert S_N(y;\tau)\vert
\le
\frac{1}{(N+1)!}
\{\frac{R(\tau+\sqrt{\tau^2+k^2})}{2}\}^{N+1}
e^{R(\tau+\sqrt{\tau^2+k^2})/2}
=E(\tau;N+1).
\tag {2.14}
$$
Now let $\tau=\tau(N)$.
Since
$$
\displaystyle
E(\tau;N)=\frac{1}{N}\frac{R(\tau+\sqrt{\tau^2+k^2})}{2}E(\tau;N-1),
$$
from (2.14) we have, as $N\longrightarrow\infty$
$$\displaystyle
\vert S_{N-1}(y;\tau(N))\vert+\vert S_{N}(y;\tau(N))\vert +\vert S_{N+1}(y;\tau(N))\vert
=O(E(\tau(N);N-1)).
\tag {2.15}
$$

\noindent
On the other hand, we have
$$
\begin{array}{c}
\displaystyle
\vert R_N(y;\tau)\vert
\le\frac{1}{(N+1)!}
\{\frac{Rk^2}{2(\tau+\sqrt{\tau^2+k^2})}\}^{N+1}
e^{\displaystyle\frac{k^2 R}{2(\tau+\sqrt{\tau^2+k^2})}}\\
\\
\displaystyle
=\frac{1}{(N+1)!}\{\frac{R(\tau+\sqrt{\tau^2+k^2})}{2}\}^{N+1}
e^{R(\tau+\sqrt{\tau^2+k^2})/2}
\{\frac{2}{R(\tau+\sqrt{\tau^2+k^2})}\}^{N+1}\\
\\
\displaystyle
\times
\{\frac{Rk^2}{2(\tau+\sqrt{\tau^2+k^2})}\}^{N+1}
e^{\{\displaystyle\frac{k^2 R}{2(\tau+\sqrt{\tau^2+k^2})}-\frac{R(\tau+\sqrt{\tau^2+k^2})}{2}\}}\\
\\
\displaystyle
=O(E(\tau;N+1))
\end{array}
$$
and thus this yields
$$\displaystyle
\vert R_{N-1}(y;\tau(N))\vert+\vert R_N(y;\tau(N))\vert+\vert R_{N+1}(y;\tau(N))\vert
=O(E(\tau(N);N-1)).
\tag {2.16}
$$

\noindent Here we claim
$$\displaystyle
e^{R\tau(N)}E(\tau(N);N-1)=O(N^{-\infty}).
\tag {2.17}
$$

\noindent
This is proved as follows.
Using the Stirling formula
$$\displaystyle
\Gamma(x)=\sqrt{\frac{2\pi}{x}}(\frac{x}{e})^x\{1+O(\frac{1}{x})\}
$$
as $x\longrightarrow\infty$ and
$$
\displaystyle
\tau(N)R+\frac{R(\tau(N)+\sqrt{\tau(N)^2+k^2})}{2}
=2\beta N/e+O(1),
$$
one gets
$$\begin{array}{c}
\displaystyle
e^{R\tau(N)}E(\tau(N);N-1)\\
\\
\displaystyle
=e^{R\tau(N)}\sqrt{\frac{N}{2\pi}}(\frac{e}{N})^{N}
\{\frac{R(\tau(N)+\sqrt{\tau(N)^2+k^2})}{2}\}^{N-1}
e^{R(\tau(N)+\sqrt{\tau(N)^2+k^2})/2}
\{1+O(\frac{1}{N})\}\\
\\
\displaystyle
=O(e^{2\beta N/e}e^{N\log(e/N)}
\sqrt{\frac{N}{2\pi}}
\{\beta N/e+O(1)\}^{N-1}\{1+O(\frac{1}{N})\})\\
\\
\displaystyle
=O(e^{\displaystyle\frac{2\beta N}{e}+N\log\frac{e}{N}
+N\log\{\beta N/e+O(1)\}})\\
\\
\displaystyle
=O(e^{\displaystyle N\{\frac{2\beta}{e}+\log\frac{e(\beta N/e+O(1))}{N}\}})
\end{array}
\tag {2.18}
$$
Since
$$\displaystyle
\frac{e(\beta N/e+O(1))}{N}
=\beta+O(\frac{1}{N}),
$$
we get
$$\displaystyle
\log\frac{\{e(\beta N/e+O(1))\}}{N}
=\log\beta+O(\frac{1}{N}).
$$
Then from (2.18) one obtains
$$\begin{array}{c}
\displaystyle
e^{R\tau(N)}E(\tau(N);N-1)
=O(e^{\displaystyle N(\frac{2\beta}{e}+\log\beta)}).
\end{array}
$$
Since
$$\displaystyle
\frac{2\beta}{e}+\log\beta<0,
$$
we obtain (2.17).

Using the recurrence relation
$$\begin{array}{c}
\displaystyle
J_{m+1}(kr)=\frac{m}{kr}J_m(kr)-J_{m}'(kr)\\
\\
\displaystyle
J_{m-1}(kr)=\frac{m}{kr}J_m(kr)+J_{m}'(kr)
\end{array}
$$
and the formulae
$$\begin{array}{c}
\displaystyle
\frac{\partial}{\partial y_1}
=\frac{e^{i\theta}}{2}(\frac{\partial}{\partial r}+i\frac{1}{r}\frac{\partial}{\partial\theta})
+\frac{e^{-i\theta}}{2}(\frac{\partial}{\partial r}-i\frac{1}{r}\frac{\partial}{\partial\theta})
\\
\\
\displaystyle
\frac{\partial}{\partial y_2}
=\frac{-ie^{i\theta}}{2}(\frac{\partial}{\partial r}+i\frac{1}{r}\frac{\partial}{\partial\theta})
+\frac{ie^{-i\theta}}{2}(\frac{\partial}{\partial r}-i\frac{1}{r}\frac{\partial}{\partial\theta}),
\end{array}
$$
one knows that: $\partial R_N(y;\tau)/\partial y_j$ can be written as a
linear combination of $R_{N+1}(y;\tau)$ and $R_{N-1}(y;\tau)$ whose coefficients
are at most algebraic growing as $\tau\longrightarrow\infty$;
$\partial S_N(y;\tau)/\partial y_j$ can be written as a linear combination
of $S_{N+1}(y;\tau)$ and $S_{N-1}(y;\tau)$ whose coefficients are at most
algebraic growing as $\tau\longrightarrow\infty$.
Using those facts, and
(2.12), (2.15), (2.16) and (2.17), one obtains the desired conclusion.

\noindent
$\Box$

\section{Proof of Theorems 1.2, 1.4 and comment on the proof of Theorem 1.5}

The starting point is the representation formula given below.
The proof is taken from that of (2.9) in \cite{CK1}.  Therein they made use of
the formula to establish an interesting equation that connects an eigenvalue of the operator
with the integral kernel $K(\varphi,d)=F(\varphi;d,k)$ acting on the functions on $S^1$
with an absorbing medium.  Here for reader's convenience we give a brief description
of the proof.

\proclaim{\noindent Lemma 3.1.} Let $\Gamma\subset S^1$ be
measurable with respect to the standard measure on $S^1$. Let
$u\in C^{\infty}(\Bbb R^2\setminus B_R)$ satisfy $\triangle
u+k^2u=0$ in $\Bbb R^2\setminus\overline B_R$; the outgoing
Sommerfeld radiation condition
$\lim_{r\longrightarrow\infty}\sqrt{r}(\partial w/\partial r
-ikw)=0$ where $r=\vert x\vert$ and $w=u-e^{ikx\cdot d}$. Then the
formula
$$\displaystyle
\int_{\Gamma}F(-\varphi;d,k)g(\varphi)d\sigma(\varphi)
=-\frac{e^{i\pi/4}}{\sqrt{8\pi k}}
\int_{\partial B_R}(\frac{\partial u}{\partial\nu}v_g-\frac{\partial v_g}{\partial\nu}u)d\sigma,
$$
is valid where $v_g$ is the Herglotz wave function with density $g\in L^2(\Gamma)$
$$\displaystyle
v_g(y)=\int_{\Gamma}e^{iky\cdot\varphi}g(\varphi)d\sigma(\varphi)
$$
and $\nu$ is the unit outward normal relative to $B_R$.
\endproclaim
{\it\noindent Proof.} From the representation formula of $w$
outside $B_R$ (\cite{CK2}) one obtains the formula
$$\displaystyle
F(\varphi;d,k)
=-\frac{e^{i\pi/4}}{\sqrt{8\pi k}}
\int_{\partial B_R}(\frac{\partial u}{\partial\nu}e^{-ik\varphi\cdot y}-\frac{\partial e^{-ik\varphi\cdot y}}
{\partial\nu}u)d\sigma(y).
$$
Thus replacing $\varphi$ with $-\varphi$, we have
$$\displaystyle
F(-\varphi;d,k)
=-\frac{e^{i\pi/4}}{\sqrt{8\pi k}}
\int_{\partial B_R}(\frac{\partial u}{\partial\nu}e^{ik\varphi\cdot y}-\frac{\partial e^{ik\varphi\cdot y}}
{\partial\nu}u)d\sigma(y).
\tag {3.1}
$$
Multiplying the both sides by $g(\varphi)$ and integrating the resultant on $\Gamma$, we obtain the desired formula.

\noindent
$\Box$

\noindent
Now we give a proof of Theorem 1.2.
Using Lemma 3.1 for $\Gamma=S^1$, we write
$$\begin{array}{c}
\displaystyle
-e^{-\tau(N) h_D(\omega)}\frac{\sqrt{8\pi k}}{e^{i\pi/4}}
\int_{S^1}F(-\varphi;d,k)g_N(\varphi;\tau(N),k,\omega)d\varphi\\
\\
\displaystyle
=e^{-\tau(N)h_D(\omega)}\int_{\partial B_R}
(\frac{\partial u}{\partial\nu}v-\frac{\partial v}{\partial\nu} u)d\sigma
+e^{-\tau(N) h_D(\omega)}\int_{\partial B_R}\{\frac{\partial u}{\partial\nu}(v_{g_N}-v)
-\frac{\partial}{\partial\nu}(v_{g_N}-v)u\}d\sigma\\
\\
\displaystyle
\equiv I_1+I_2
\end{array}
\tag {3.2}
$$
where $v=e^{y\cdot(\tau\omega+i\sqrt{\tau^2+k^2}\omega^{\perp})}$ and
$v_{g_N}$ is the Herglotz wave function with density $g_N(\,\cdot\,;\tau(N),k,\omega)$.

\noindent
In \cite{I3} we have already proven that there exist $\mu>0$ and
$A>0$ such that
$$\displaystyle
\lim_{N\longrightarrow\infty}\tau(N)^{\mu}\vert I_1\vert=A.
\tag {3.3}
$$
Theorem 2.1 gives the estimate
$$\begin{array}{c}
\displaystyle
\tau(N)^{\mu}\vert I_2\vert\\
\\
\displaystyle
\le \tau(N)^{\mu}e^{R\tau(N)}C_R(u)\{\sup_{\vert y\vert\le R}\vert v_{g_N}(y)-v(y)\vert
+\sup_{\vert y\vert\le R}\vert\nabla\{v_{g_N}(y)-v(y)\}\vert\}
=O(N^{-\infty})
\end{array}
\tag {3.4}
$$
where $C_R(u)$ is a positive constant depending on $u$, $R$ and independent of $N$.

\noindent
From (3.2), (3.3) and (3.4) one gets
$$\displaystyle
\tau(N)^{\mu}e^{-\tau(N) h_D(\omega)}
\vert\int_{S^1}F(-\varphi;d,k)g_N(\varphi;\tau(N),k,\omega)d\sigma(\varphi)\vert\longrightarrow\frac{A}{\sqrt{8\pi k}}
\tag {3.5}
$$
as $N\longrightarrow\infty$.
Then (1.4) and other all conclusions come from (3.5).

\noindent
$\Box$

In the following we say that: a function $f(\tau)$ decays algebraically as $\tau\longrightarrow$
in strict sense if $\tau^{\lambda}\vert f(\tau)\vert$ converges to a positive number
as $\tau\longrightarrow$ for a positive constant $\lambda$;
$f(\tau)$ is decaying at most algebraic if $f(\tau)=O(\tau^{-\mu})$ as $\tau\longrightarrow\infty$
for a suitable positive constant $\mu$.

\noindent
Now we give a proof of Theorem 1.4.
Let $y$ denote the only one point of the set
$\{x\,\vert\,x\cdot\omega=h_{\Sigma}(\omega)\}\cap\Sigma$. First
consider the case when every end points of
$\Sigma_1,\Sigma_2,\cdots,\Sigma_m$ satisfies
$x\cdot\omega<h_{\Sigma}(\omega)$. Then $y$ should be a point
where two segments in some $\Sigma_j$ meet.  Then from Theorem 4.2
in \cite{I3} and a fact similar to Lemma 5.1 in \cite{I3} one obtains that both
$I^{d_1}_{\omega}(\tau,h_{\Sigma}(\omega))$ and
$I^{d_2}_{\omega}(\tau,h_{\Sigma}(\omega))$ decays algebraically
as $\tau\longrightarrow\infty$ in strict sense.  This yields the
algebraic decaying of $I_{\omega}(\tau, h_{\Sigma}(\omega))$ as
$\tau\longrightarrow\infty$ in strict sense. Then we automatically
obtain the desired results.

The problem is the case when $y$ is an end point of some
$\Sigma_j$. Since $\nu$ on $\Sigma_j$ becomes a constant vector in
a neighbourhood of $y$, we denote the constant vector by
$\nu_j$.  Using Theorem 4.3 in \cite{I3} for $u$ with $d=d_1, d_2$
one concludes that $I_{\omega}(\tau,h_{\Sigma}(\omega))$ decays at
most algebraically as $\tau\longrightarrow\infty$. Moreover, an
argument similar to that of the proof of Lemma 5.1 in \cite{I3}
yields that if both of $I^{d_1}(\tau, h_{\Sigma}(\omega))$ and
$I^{d_2}(\tau, h_{\Sigma}(\omega))$ decay rapidly then
$d_1\cdot\nu_j=d_2\cdot\nu_j=0$.   However, this is impossible.
Thus one concludes that one of them has to be algebraically
decaying as $\tau\longrightarrow\infty$ in strict sense.  Then we
obtain the algebraic decaying of $I_{\omega}(\tau,h_{\Sigma}(\omega))$
as $\tau\longrightarrow\infty$ in strict sense.
This completes the proof.

\noindent
$\Box$

\noindent
It is easy to see that completely the same formula for $d=d_1, d_2$ as Lemma 3.1 is valid.
Then the proof of Theorem 1.5 can be done along the same line as that of Theorem 1.2.

\section{Limited aperture}

In this section we consider the case when the far field pattern with limited aperture is given.
Here we point out an effect of a priori information
on the formulae in Theorem 1.2.
Let $\Gamma $ be a non empty open subset of $S^1$.

\noindent
We assume that:

(1) the far field pattern on $\Gamma$ is known for fixed $d$ and $k$;

(2) $0\in D$.

\noindent
(1) means that $\Gamma$ is an aperture.  (2) means that the center of the coordinates is inside $D$
and we know it in advance.

The starting point is a theorem established in \cite{CK3}.
Define
$$\displaystyle
W(B_R)=\{v\in C^2(B_R)\cap C^1(\overline B_R)\,\vert\,\triangle v+k^2v=0\,\,\text{in}\,B_R\}.
$$
We denote by $\overline{W(B_R)}$ the $H^1(B_R)$ closure of $W(B_R)$.
Given $g\in L^2(-\Gamma)$ define
$$\displaystyle
Hg(y)=\int_{-\Gamma}e^{iky\cdot\varphi}g(\varphi)d\sigma(\varphi)\,\,\,y\in\,B_R.
$$
Then $Hg\in \overline{W(B_R)}$ and the map $H:L^2(-\Gamma)\longrightarrow\overline{W(B_R)}$
is bounded linear.  They proved the following.

\proclaim{\noindent Theorem 4.1(Theorem 2.6 of \cite{CK3}).}
The range of $H$ is dense in $\overline{W(B_R)}$.
\endproclaim

\noindent
Note that: therein they considered only the case when $-\Gamma=S^1$, however,
by using the real analyticity of the far field pattern, one knows that the proof is still valid.

Given $v\in \overline{W(B_R)}$ and $\delta>0$ an element $g_0\in L^2(-\Gamma)$ is called
minimum norm solution of $Hg=v$ with discrepancy $\delta$, if $g_0$ satisfies
$\Vert Hg-v\Vert_{H^1(B_R)}\le\delta$ and
$$\displaystyle
\Vert g_0\Vert_{L^2(-\Gamma)}
=\inf\{\Vert g\Vert_{L^2(-\Gamma)}\,:\,
\Vert Hg-v\Vert_{H^1(B_R)}\le\delta\}.
$$
Now given $\tau>0$ and $\omega\in S^1$ set $v=e^{x\cdot(\tau\omega+i\sqrt{\tau^2+k^2}\omega^{\perp})}$($\in\overline {W(B_R)}$).
Theorem 4.1 ensures the existence of the minimum norm solution of $Hg=v$ with discrepancy $\delta$ (see Theorem 16.11 in
\cite{Kr}).  It is given by the formula
$$\displaystyle
g=(\alpha I+H^*H)^{-1}H^*v
$$
where $\alpha>0$ is any zero of the function
$$
\Vert H(\alpha I+H^*H)^{-1}H^*v-v\Vert_{H^1(B_R)}^2-\delta^2.
$$
We denote the minimum norm solution by $g=g_{\tau,\delta}(\,\cdot\,;k,\omega)$.
This satisfies
$$\displaystyle
\Vert Hg_{\tau,\delta}(\,\cdot\,;k,\omega)-v\Vert_{H^1(B_R)}\le\delta.
\tag {4.1}
$$
Then we obtain the following theorem.

\proclaim{\noindent Theorem 4.2.}
Assume that $0\in D$.
Let $\omega$ be regular with
respect to $D$.  Then the formula
$$\displaystyle
\lim_{\tau\longrightarrow\infty}
\frac{\displaystyle\log\vert\int_{\Gamma}F(\varphi;d,k)g_{\tau,\delta}(-\varphi;k,\omega)d\sigma(\varphi)\vert}
{\tau}=h_D(\omega),
$$
is valid.
Moreover we have:

if $t\ge h_D(\omega)$, then
$$\displaystyle
\lim_{\tau\longrightarrow\infty}e^{-\tau t}
\vert\int_{\Gamma}F(\varphi;d,k)g_{\tau,\delta}(-\varphi;k,\omega)d\sigma(\varphi)\vert=0;
$$

if $t<h_D(\omega)$, then
$$\displaystyle
\lim_{\tau\longrightarrow\infty}e^{-\tau t}
\vert\int_{\Gamma}F(\varphi;d,k)g_{\tau,\delta}(-\varphi;k,\omega)d\sigma(\varphi)\vert=\infty.
$$

\endproclaim

{\it\noindent Proof.}
Using Lemma 3.1, we write
$$\begin{array}{c}
\displaystyle
-e^{-\tau h_D(\omega)}\frac{\sqrt{8\pi k}}{e^{i\pi/4}}
\int_{\Gamma}F(\varphi;d,k)g_{\tau,\delta}(-\varphi;k,\omega)d\varphi\\
\\
\displaystyle
=e^{-\tau h_D(\omega)}\int_{\partial B_R}
(\frac{\partial u}{\partial\nu}v-\frac{\partial v}{\partial\nu} u)d\sigma
+e^{-\tau h_D(\omega)}\int_{\partial B_R}\{\frac{\partial u}{\partial\nu}(Hg_{\tau,\delta}-v)
-\frac{\partial}{\partial\nu}(H{g_{\tau,\delta}}-v)u\}d\sigma
\end{array}
\tag {4.2}
$$
where $v=e^{y\cdot(\tau\omega+i\sqrt{\tau^2+k^2}\omega^{\perp})}$ and
$Hg_{\tau,\delta}$ is the Herglotz wave function with density $g_{\tau,\delta}(\,\cdot\,;k,\omega)$.
Using the trace theorem, from (4.1) we know that
$$
\displaystyle
\Vert Hg_{\tau,\delta}-v\Vert_{H^{1/2}(\partial B_R)}
+\Vert \frac{\partial}{\partial\nu}
\{Hg_{\tau,\delta}-v\}\Vert_{H^{-1/2}(\partial B_R)}\le C\delta
\tag {4.3}
$$
where $C>0$ is independent of $\tau$.  The assumption gives
$h_D(\omega)>0$. Thus from (4.3) one knows that the second term of
the right hand side of (4.2) is exponentially decaying as
$\tau\longrightarrow\infty$ and of course, the first term is
decaying algebraically as $\tau\longrightarrow\infty$ in strict
sense by the same reason as described in the proof of Theorem 1.2.

\noindent
$$\Box$$

\noindent
In a forthcoming paper we consider
algorithms that are based on Theorems 1.2, 1.5 and 4.2.

$$\quad$$

\centerline{{\bf Acknowledgement}}

This research was partially supported by Grant-in-Aid for Scientific
Research (C)(2) (No. 15540154) of Japan Society for the Promotion
of Science.  The author thanks the anonymous referee for suggestions for improvement
of the original manuscript.

\vskip1cm
\noindent
e-mail address

ikehata@math.sci.gunma-u.ac.jp
\end{document}